\documentclass[11pt]{article}
\usepackage[usenames]{color}
\usepackage{epsf}
\usepackage{amsmath, amsfonts, amssymb, amsthm}
\usepackage{fullpage}
\usepackage{hyperref}

\newtheorem{theorem}{Theorem}[section]

\newtheorem{lemma}[theorem]{Lemma}

\newtheorem{claim}[theorem]{Claim}
\newtheorem{example}[theorem]{Example}

\title{Balanced Spanning Caterpillars}
\author{Andrzej Czygrinow\thanks{Partially supported by Simons Foundation Grant \# 521777} and Jangwon Yie\\School of Mathematical and Statistical Sciences\\Arizona State University\\Tempe, AZ 85287}
\begin{document}
\maketitle

\begin{abstract}
A $p$-caterpillar is a caterpillar such that every non-leaf vertex is adjacent to exactly $p$ leaves. We give a tight minimum degree condition for a graph to have a spanning $p$-caterpillar.
\end{abstract}
\section{Introduction}
Every connected graph contains a spanning tree, yet quite often it is desirable to find a spanning tree which satisfies certain additional conditions. 
There are many results giving sufficient minimum degree conditions for graphs to contain very special spanning trees. For example, Dirac's theorem  from \cite{dirac} states that any graph on $n\geq 3$ vertices with minimum degree at least $(n-1)/2$ has a spanning path. In \cite{win}, S. Win  generalized this fact and proved the following theorem.
\begin{theorem}\label{win-thm}
Let $k\geq 2$ and let $G$ be a graph on $n$ vertices such that $\sum_{x\in I} d(x)\geq n-1$ for every independent set $I$ of size $k$. Then $G$ contains a spanning tree of maximum degree at most $k$. 
\end{theorem}
In particular, if the minimum degree of $G$ is at least $(n-1)/k$, then $G$ contains a spanning tree of maximum degree at most $k$. In fact, as showed in \cite{czyg}, the degree condition from Theorem \ref{win-thm} implies that either $G$ has a spanning caterpillar of maximum degree at most $k$ or $G$ belongs to a special exceptional class. We refer the reader to \cite{ozeki} for a comprehensive survey of spanning trees. 

Another way of thinking about caterpillars is by looking at domination problems. A set $S\subseteq V$ is a dominating set in a graph $G=(V,E)$ if every vertex in $V\setminus S$ has a neighbor in $S$. A dominating set $S$ is called a {\it connected dominating set} if, in addition, $G[S]$ is connected. 
In the special case when $G[S]$ contains a path, we say that $G$ has a dominating path. In \cite{bro}, Broersma proved a result on cycles passing within a specified distance of a vertex and stated an analogous result for paths from which, as one of the corollaries, we get the following fact.
\begin{theorem}
If $G$ is a $k$-connected graph on $n$ vertices such that $\delta(G)> \frac{n-k}{k+2}-1$, then $G$ contains a dominating path.
\end{theorem}
In particular, if $G$ is connected then $\delta(G)>\frac{n-1}{3}-1$ implies that $G$ has a spanning caterpillar.
In this paper we will be concerned with a minimum degree condition that implies existence of spanning balanced caterpillar. 

A $p$-caterpillar is a tree such that the graph induced by its internal vertices is a  path and every internal vertex has exactly $p$ leaves. The {\it spine} of a caterpillar is the graph induced by its internal vertices. The length of a caterpillar is the length of its spine.
Recently, Faudree et. al. proved the following fact in \cite{faudree}.
\begin{theorem}\label{west}
For $p\in Z^+$ there exists $n_0$ such that for every $n\in (p+1)Z$ such that $n\geq n_0$ the following holds. If $G$ is a graph on $n$ vertices such that $\delta(G)\geq \left(1-\frac{p}{(p+1)^2}\right) n$, then $G$ contains a spanning $p$-caterpillar. 
\end{theorem}
The authors of \cite{faudree} ask for the tight minimum degree condition which implies that $G$ has a spanning $1$-caterpillar. In addition, they ask for a tight minimum degree condition which gives a nearly balanced $p$-caterpillar (every vertex on the spine has $p$ or $p+1$ leaf neighbors).
We will settle the first problem and answer the second question in the case when $n$ is divisible by $p+1$.
In this paper we will substantially improve the minimum degree bound from Theorem \ref{west} and give a tight minimum degree condition which guarantees existence of a spanning $p-$caterpillar.
Our main result is the following fact.
\begin{theorem}\label{main-theorem}
For $p\in Z^+$, there exists $n_0$ such that for every $n \in (p+1)Z$ with $n \geq n_0$ the following holds. If $G$ is a graph on $n$ vertices such that
$$\delta(G)\geq \left\lbrace 
\begin{array}{ll} \frac{n}{2}   & \text{ if } n/(p+1) \text{ is even} \\
  \frac {n+1}{2} & \text{ if }  n/(p+1) \text{ is odd and } p> 2\\
   \frac{n-1}{2} & \text { if } n/(p+1) \text{ is odd and } p\leq 2 \end{array}\right.
$$
then $G$ contains a spanning $p$-caterpillar.
\end{theorem}
It's not difficult to see that the minimum degree condition in Theorem \ref{main-theorem} is best possible.
\begin{example}
First note that $K_{n/2}\cup K_{n/2}$ in the case $n$ is even and $K_{(n-1)/2} \cup K_{(n+1)/2}$ in the case $n$ is odd have no spanning caterpillars. Thus the degree condition in the case $p\leq 2$ is tight.
Now suppose $p\geq 3$. Let $n/(p+1)$ be even. Then $n/2$ is an integer. Consider $K_{n/2-1, n/2+1}$. Clearly  $n/(2(p+1))$ of spine vertices must be in one of the partite sets, because the spine is a path and its maximum independent set is of size $n/(2(p+1))$, but then the two partite sets must have the same size. Another example is $K_{n/2}\cup K_{n/2}$.
Now, suppose $n/(p+1)$ equals $2k+1$ for some $k\in Z^+$. If $n$ is even, then consider $K_{n/2, n/2}$. Clearly one of the partite sets must have $k+1$ spine vertices and so the other set must contain $(k+1)(p+1) - 1 = \frac{n+p-1}{2}> n/2$ as $p>1$. If $n$ is odd then consider $K_{(n-1)/2, (n+1)/2}$. 
Now, $k+1$ of the spine vertices must be in the partite set of size $(n-1)/2$. Consequently, the other set must have at least $\frac{n+p-1}{2} > (n+1)/2$ as $p>2$. 
\end{example}

We will prove Theorem \ref{main-theorem} using the absorbing method from \cite{RRS3}. In this method, we first analyze the non-extremal case and then address two extremal cases, when  $G$ is "close to" $2K_{\lfloor n/2\rfloor}$ or $K_{ \lfloor n/2\rfloor, \lceil n/2 \rceil}$.

We will use $|G|, ||G||$ to denote the order and the size of a graph $G$. For two, not necessarily disjoint, sets $U,W\subseteq V(G)$, we will use $||U,W||$ to denote the number of edges in $G$ with one endpoint in $U$, another in $W$.
We say that a graph $G$ is $\beta$-extremal if either $V(G)$ contains a set $W$ such that $|W|\geq (1/2-\beta)n$ and $||G[W]||\leq \beta n^2$ or if $V(G)$ can be partitioned into sets $V_1,V_2$ so that $|V_i|\geq (1/2-\beta) n$ for  $i=1,2$ and $||V_1,V_2||\leq \beta n^2$. 
In addition, the following notation and terminology will be used. A $u,v$-caterpillar is a $p$-caterpillar where the first vertex in the spine is $u$ and the last is $v$. 

The rest of the paper is structured as follows. In Section \ref{abs-section} we prove the absorbing lemma which is the key to handle the non-extremal case. In Section \ref{non-extremal-sec} we prove the non-extremal case and in Section \ref{extremal-sec} we address the extremal cases.

\section{Absorbing Lemma}\label{abs-section}
In this section we will prove an absorbing lemma and a few additional facts which are used in the next section to complete the proof in the case a graph is not extremal. We will start with the following observation.
\begin{lemma}\label{no-edges}
For $1/8>\beta>0$ there is $\alpha >0$ and $n_0$ such that the following holds. If $G$ is a graph on $n \geq n_0$ vertices such that $\delta(G)\geq (1/2-\beta^2)n$ which is not $\beta$-extremal, then for any (not necessarily distinct) vertices $u,v\in G$, $||N(u), N(v)|| \geq \beta^2n^2/32$.
\end{lemma}
{\bf Proof.} We have $||G[N(u)]|| > \beta n^2$ from the definition of a $\beta$-extremal graph. Now suppose $u,v$ are two distinct vertices. If $\beta n/2 \leq |N(u)\cap N(v)| \leq (1/2 -\beta/2)n$, then $|N(u) \cup N(v)| \geq 2(1/2-\beta^2)n -(1/2-\beta/2)n \geq (1/2+ \beta/4)n$. Thus every vertex $x\in N(u)\cap N(v)$ has at least $\beta n/8$ neighbors in $N(u)\cup N(v)$. Consequently, $||N(v), N(u)|| \geq \beta^2 n^2/32$. If $|N(u)\cap N(v)|< \beta n/2$, then $|N(v)\setminus N(u)| \geq (1/2 - 2\beta/3)n$. Thus, since $G$ is not $\beta$-extremal $||N(u), N(v)||\geq \beta^2 n^2/32$. If $|N(u) \cap N(v)| \geq (1/2 -\beta/2)n$, then  $||G[N(u) \cap N(v)]|| \geq \beta n^2$. $\Box$

Our next objective is to establish the following connecting lemma.
\begin{lemma}[Connecting Lemma]\label{connect}
For $1/8 > \beta >0$ there is $\alpha >0$ and $n_0$ such that the following holds. If $G$ is a graph on $n \geq n_0$ vertices such that $\delta(G)\geq (1/2-\beta^2)n$ which is not $\beta$-extremal, then for any two vertices $u,v\in G$ there are at least $\alpha n^{4p+2}$ $u,v$-caterpillars of length three in $G$.
\end{lemma}
{\bf Proof.} Let $u,v$ be two distinct vertices. By Lemma \ref{no-edges}, $||N(u), N(v)||\geq \beta^2n^2/32$. Let $\{x,y\} \in E(N(u), N(v))$. Since each vertex in $\{x,y,u,v\}$ has degree at least $(1/2-\beta^2)n$, the number of different $p$-caterpillars with spine $u,x,y,v$ is at least $\gamma n^{4p}$ for some $\gamma>0$. Thus the total number of $u,v$-caterpillars of length three in $G$ is at least $\alpha n^{4p+2}$ for some $\alpha>0$ which depends on $\beta$ only. $\Box$

We will be connecting through a small subset of $V(G)$ called a reservoir set.
\begin{lemma}[Reservoir Set]\label{reservoir-lem}
For $1/64 > \beta >0$ and $\beta^4> \gamma >0$ there is $n_0$ such that if $G$ is a graph on $n \geq n_0$ vertices satisfying $\delta(G)\geq (1/2-\beta^2)n$ which is not $\beta$-extremal then there is a set $Z\subset V(G)$ such that  the following holds:
\begin{itemize}
\item[(i)] $|Z| = (\gamma \pm \gamma^2) n$;
\item[(ii)] For every $v\in V$, $|N(v) \cap Z| \geq (1/2- 2\beta^2)\gamma n$;
\item[(iii)] For every $u,v\in V$, $||N(u) \cap Z, N(v) \cap Z|| \geq \beta^6 \gamma^2 n^2/4$.
\end{itemize}
\end{lemma}
{\bf Proof.} Let $Z$ be a set obtained by selecting every vertex from $V$ independently with probability $p := \gamma$.
By the Chernoff bound \cite{chernoff}, with probability $1-o(1)$, the following facts hold:
\begin{itemize}
\item[(a)] $(\gamma - \gamma^2)n \leq |Z| \leq (\gamma + \gamma^2)n$ ;
\item[(b)] For every vertex $v$, $|N(v) \cap Z| \geq (1/2 - 2\beta^2)\gamma n$.
\end{itemize}

To prove the third part let $u,v\in V$ and let $X_{u,v}:=\{w\in N(u)| |N(w) \cap N(v)|\geq \beta^3 n\}$.
Since $G$ is not $\beta$-extremal by Lemma \ref{no-edges}, $||N(u), N(v)||\geq \beta^2n^2/32$.
Thus $|X_{u,v}|\geq \beta^3 n$. Indeed, if $|X_{u,v}| < \beta^3n $, then $||N(u),N(v)|| < 2\beta^3n^2 <  \beta^2n^2/32$. Consequently, by Chernoff's inequality, with probability $1- o(1/n^2), |X_{u,v} \cap Z| \geq \beta^3 \gamma n/2$. Thus with probability $1-o(1)$ for every $u,v$, $|X_{u,v}|\geq \beta^3\gamma n/2$.
Let $u\in V$ be arbitrary  and let  $w \in V$ be such that $|N(w) \cap N(u)| \geq \beta^3 n$. Then with probability at least $1 - o(1/n^2)$, $|N(w) \cap N(u) \cap Z|\geq \beta^3 \gamma n/2$.
Thus with probability at least $1 - o(1)$, we have
$$||N(u) \cap Z, N(v) \cap Z|| \geq \beta^6 \gamma ^2 n^2/4$$
for every $u,v$.
Therefore there is a set $Z$ such that (i)-(iii) hold. $\Box$

We will continue with our proof of the absorbing lemma.
We shall assume that $0 < \beta< 1/64$, $G=(V,E)$ is a graph on $n$ vertices where $n$ is sufficiently large which is not $\beta$-extremal and which satisfies $\delta(G) \geq (n-1)/2$.
In addition, we will use an auxiliary constant $\tau$ such that $0< \tau < \frac{\beta}{10}$.
\begin{lemma}\label{lemma-cond} Let $u,v$ be two vertices in $G$ such that $|N(u) \cap N(v)|\geq 2\tau n$. Then, at least one of the following conditions holds.
\begin{itemize}
\item[(1)] At least $\tau n$ vertices $x\in N(u) \cap N(v)$ are such that $|N(x)\cap N(u)|\geq \tau^2 n$.
\item[(2)] All but at most $3\tau n$ vertices $x\in N(v)$ satisfy $|N(x)\cap N(v)|\geq \tau^3 n$.
\end{itemize}
\end{lemma}
{\bf Proof.} First suppose that $|N(v) \setminus N(u)|< 2\tau n + 2$. Since $G$ is not $\beta$-extremal, $|| G[N(v) \cap N(u)]|| \geq \beta n^2$ and so the first condition holds. Thus we may assume that $|N(v) \setminus N(u)|\geq 2\tau n + 2$.
Since $|N(u)\cup N(v)| > (1/2+2\tau) n + 1$, every vertex $x \in N(v)\cap N(u)$ has at least $2\tau n$ neighbors in $N(u) \cup N(v)$.
Thus all but at most $\tau n$ vertices in $N(v) \cap N(u)$ have at least $(2\tau - \tau^2) > \tau^3n$ neighbors in $N(v)$.

Now, suppose the first condition fails and we claim that all but at most $2\tau n$ vertices $x \in N(v) \setminus N(u)$ satisfy $|N(x) \cap N(v)| \geq \tau^3 n$.
Let $A:=\{x\in N(u) \cap N(v) : |N(x) \cap N(u)| < \tau^2 n\}$ and note that $|A| \geq \tau n$.
Therefore,
$$||A, N(v) \setminus N(u)|| \geq |A|(|N(v)\setminus N(u)|-\tau^2n) \geq (1-\tau)|A||N(v)\setminus N(u)|.$$
Let $B:=\{y \in N(v)\setminus N(u) : |N(y)\cap A|< \tau^2 |A|\}$ then $\tau^2|B||A| +(|N(v)\setminus N(u)|-|B|)|A| > ||A, N(v) \setminus N(u)||$.
Hence $\tau^2|B||A| +(|N(v)\setminus N(u)|-|B|)|A| > (1-\tau)|A||N(v)\setminus N(u)|$ and so $|B| < \tau |N(v)\setminus N(u)|/(1-\tau^2) < 2\tau n$.
For every vertex $x\in (N(v)\setminus N(u))\setminus B$,
$$|N(x) \cap N(v)|\geq |N(x)\cap A|\geq \tau^2 |A| \geq \tau^3 n,$$
which completes the proof.$\Box$
\begin{lemma}\label{lemma-exists}
Let $T$ be a set of $p+1$ vertices in $G$. Then there exists a vertex $x\in T$ such that for every $y\in T$, $||N(x), N(y)||\geq \tau^4 n^2$.
\end{lemma}
{\bf Proof.} Suppose there is a vertex $v\in T$ such that condition (2) in Lemma \ref{lemma-cond} is satisfied. Let $x:=v$ and take $y\in T$.
If $|N(y)\cap N(x)|\geq 5\tau n$, then $||N(x), N(y)||\geq \frac{1}{2} \cdot (2\tau^4 n^2)$. If $|N(y) \cap N(x)| < 5\tau n$, then since $G$ is not $\beta$-extremal, $||N(x), N(y)|| \geq \tau^4 n^2$.
Therefore, we may assume that there is no such $v$ in $T$. Let $x$ be an arbitrary vertex in $T$. Take $y\in T$. If $|N(x) \cap N(y)|\geq 2\tau n$, then by Lemma \ref{lemma-cond} (with $u:=x$, $v:=y$), $||N(x), N(y)||\geq \tau^3 n^2/2$. If $|N(x)\cap N(y)|< 2\gamma n$, then since $G$ is not $\beta$-extremal, $||N(x), N(y)||\geq \gamma^4 n^2$. $\Box$

We say that an $x,y$-caterpillar $P$ absorbs a set $T$ of size $p+1$, if $G[V(P)\cup T]$ contains an $x,y$-caterpillar on $|V(P)|+p+1$ vertices. Let $M_{q}(T)$  denote the set of caterpillars of order $q$ which absorb $T$.
A caterpillar $P$ is called {\it $\gamma-$absorbing} if $P$ absorbs every subset $W \subset V \setminus V(P)$ with $|W|\in (p+1)Z$ and $|W| \leq \gamma n$.
We will now prove our main lemma from which the absorbing lemma follows by using the deletion method.
 \begin{lemma}\label{lem-count}
Let $p \in Z^+$. For every $\beta>0$ there is $n_0$ and $\alpha>0$ such that the following holds. If $G$ is a graph on $n \geq n_0$ vertices which is not $\beta$-extremal and such that $\delta(G)\geq (n-1)/2$ and $T \subset V(G)$, $|T|=p+1$, then
 $$ |M_{q}(T)| \geq \alpha n^q$$
 where $q=(3p+2)(p+1)$.
 \end{lemma}
 {\bf Proof.} Let $T=\{x,y_1, \dots, y_p\}$ and in view of Lemma \ref{lemma-exists} suppose that for every $i$, $||N(x), N(y_i)||\geq \tau^4 n^2$.
  We will construct a caterpillar $P$ which absorbs $T$. The counting fact follows easily from the way the construction works.
To construct the caterpillar we will proceed in a few steps, selecting distinct vertices which have not been previously selected in each step.
First take $v_i \in N(y_i)$ so that $v_1, \dots, v_p$ are distinct and $|N(v_i) \cap N(x)|\geq \tau^5 n$.
Now let $u_i \in N(v_i) \cap N(x)$  be such that $u_1, \dots, u_p$ are distinct. Let $x_1x_2$ be an edge in $N(x)$.
Use Lemma \ref{connect}, to find $v_i,v_{i+1}$ caterpillars with spines $P_i$ for $2 \leq i \leq p-1$, all vertices distinct, and let $P:=v_2P_2v_3\dots v_{p-1}P_{p-1}v_p$.
Use Lemma \ref{connect} to find a $v_2, x_2$-caterpillar and denote its spine by $Q_2$ and a $v_1, x_1$-caterpillar with spine $Q_1$. Let $Q:= v_1Q_1x_1x_2Q_2v_2Pv_p$.
Then $Q$ is a $v_1,v_p$-path.  Disregard selected vertices not on $Q$. For every vertex $v_i$ select $p-1$ distinct neighbors, so that together with $u_i$ they give $p$ leaves attached to $v_i$. For $x_1, x_2$ select $p$ distinct neighbors and let $S$ be the set containing all the vertices on $Q$,  $u_1, \dots, u_p$, and all the remaining neighbors. Then $G[S]$ contains a $v_1,v_p$- caterpillar of length $3p+1$ which contains  $(3p+2)(p+1)$ vertices.
In addition, $G[S \cup T]$ contains a $v_1,v_p$-caterpillar of length $3p+2$ obtained as follows. Insert $x$ between $x_1$ and $x_2$ in the spine $Q$, make $u_1, \dots, u_p$ the neighbors of $x$, and let $y_i$  replace $u_i$ in the set of spikes of $v_i$. By Lemma \ref{connect} and in view of the construction, the number of such sets $S$ is at least $\alpha n^{(3p+2)(p+1)}$ for some $\alpha>0$ which depends on $\beta$ and $p$ only. $\Box$

\begin{lemma}\label{absorbing-lemma}
(Absorbing Lemma) Let $p \in Z^+, q = (3p+2)(p+1), \beta > 0$ and $\alpha > 0 $ be such that Lemma \ref{lem-count} holds. For any $\delta < \alpha / 10q$, there is $n_0$ such that the following holds.
If $G$ is a graph on $n \geq n_0$ vertices which is not $\beta$-extremal and such that $\delta(G)\geq (n-1)/2$ then there is a caterpillar $P_{abs}$ in $G$ on at most $\delta n$ vertices which is $\delta^2$-absorbing.
\end{lemma}
{\bf Proof.}
Let $n_0$ be such that Lemma \ref{lem-count} holds with $\alpha$.
Let $G$ be a graph on $n \geq n_0$ vertices which is not $\beta-$extremal and such that $\delta(G) \geq (n-1)/2$.

Let $\mathcal{F}$ be a family obtained by selecting every set from $\binom{V}{q}$ independently with probability $\mu := \delta n / 3q \binom{n}{q}$.
By the Chernoff bound \cite{chernoff}, with probability $1 - o(1)$,
$$| \mathcal{F}| \leq 2\mu \binom{n}{q} = 2\delta n / 3q$$
Now, let $T$ be a set of size $p+1$. Again by Chernoff bound, with probability $1 - o(1/n^{p+1}),$
$$ |M_q(T) \cap \mathcal{F}| \geq \frac{1}{2} \mu \alpha n^q > 3\delta^2 n. $$

The expected number of pairs $\{S_1,S_2\}$ such that $S_1,S_2 \in \mathcal{F}$ and $S_1 \cap S_2 \neq \emptyset$ is at most $p\binom{n}{q} \cdot q \binom{n}{q-1}p \leq \delta^2 n$ and so by Markov's inequality, with probability at least 1/2, the number of such pairs is at most $2 \delta^2 n$.
Therefore, with positive probability, there exists a family $\mathcal{F}$ such that $|\mathcal{F}| \leq 2\delta n / 3q$, for every set $T$ of size $p+1$, $|M_q(T) \cap \mathcal{F}| > 3 \delta^2 n$, and the number of $\{S_1,S_2\}$ such that $S_1,S_2 \in \mathcal{F}$ and $S_1 \cap S_2 \neq \emptyset$ is at most $2 \delta^2 n$.
Let $\mathcal{F}'$ be obtained from $\mathcal{F}$ by deleting all intersecting sets and sets that do not absorb any $T$. Then $|\mathcal{F}'| \leq 2 \delta n /3q$, and for every set $T$ of size $p+1$, $|M_q(T) \cap \mathcal{F}'| > \delta^2 n$. For each $S \in \mathcal{F}'$, $G[S]$ contains a caterpillar on $q$ vertices, so by using the minimum degree condition and Lemma \ref{connect}, we can connect the endpoints of these caterpillars to obtain a new caterpillar $P_{abs}.$ We also have that
$$ |P_{abs}| \leq |\mathcal{F}'| \cdot q + 2|\mathcal{F}'| \cdot p < |\mathcal{F}'| \cdot (3q /2) \leq \delta n.$$

To show that $P_{abs}$ is $\delta^2-$absorbing, consider $W \subset V \setminus V(P_{abs})$ such that $(p+1) | |W|$ and $|W| \leq \delta^2 n$. $\mathcal{W} = \{ W_1, ... ,W_m \}$ be an arbitrary partition of $W$ into sets of size $p+1$. We have that $|M_q(W_i) \cap \mathcal{F}'| > \delta^2 n $ for every $i \in [m]$. Therefore, there exists a matching between $\mathcal{W}$ and $\mathcal{F}'$ so that every $W_i \in \mathcal{W}$ is paired with some $S_i \in M_q(W_i)$. This implies that $P_{abs}$ absorbs $W$ and the proof is complete.
$\Box$
\section{Non-extremal case}\label{non-extremal-sec}
In this section we will finish proving the non-extremal case. The argument uses a similar approach as the proof of a corresponding fact in \cite{CM}.

Let $p \in Z^{+}, q = (3p+2)(p+1)$ and let $\xi, \beta $ be such that $0 < \xi < 1/(4p+5), 0 < \beta < \min \{ (\frac{\xi}{30p})^2, (\frac{\xi}{96})^2\}$.
Now, let $\alpha >0, n_0 \in \mathbb{N}$ be such that Lemma \ref{absorbing-lemma} holds.
Let $\delta, \gamma > 0$ be such that $\delta < \min \{ (\frac{\beta}{300})^2, \frac{\alpha}{10q} \}, \gamma < \frac{\delta^2}{4}$ and $C$ be such that $C > \frac{80(p+1)}{\delta \gamma \beta^3}$.
Let $n > \max \{n_0, \frac{4C \cdot 2^{(1+\delta)C} }{\delta^3} \}$ and $G$ be a graph on $n$ vertices which is not $\beta-$extremal and of minimum degree at least $(n-1)/2$.
Let $P_{abs}$ be the absorbing caterpillar obtained in the previous section and let $Z$ be the reservoir set from Lemma \ref{reservoir-lem} applied with $\gamma$ which is less than $\beta^4$ because  $\gamma < \frac{\delta^2}{4} <\beta^4 $.

\begin{claim} \label{gluing}
Let $P_1,P_2$ be disjoint caterpillars in $G$ such that $|Z \cap V(P_1)|, |Z \cap V(P_2)| < \frac{\beta^3 \gamma n}{4}$ and the  endpoints of $P_1$ and $P_2$ are not in $Z$.
Then there is a caterpillar $P$ containing $V(P_1)\cup V(P_2)$ which has at most $2(p+1)$ additional vertices in $Z$ and such that its endpoints are not in $Z$.
\end{claim}
{\bf Proof.}
Let $u_1,u_2$ be the endpoints of $P_1,P_2,$ respectively.
Since $\frac{\beta^3}{2} < 1/2 -2 \beta^2$ and $|(N(u_1) \cap Z) \cap (V(P_1) \cup V(P_2))| \cdot |(N(u_2) \cap Z) \cap (V(P_1) \cup V(P_2))| < \beta^6 \gamma^2 n^2 / 16$, by Lemma \ref{reservoir-lem}, there exists $ x_1 \in (N(u_1) \cap Z) - (V(P_1) \cup V(P_2)) , x_2 \in (N(u_2) \cap Z) - (V(P_1) \cup V(P_2))$ such that $\{x_1, x_2\} \in E(G)$.
Then we can construct new caterpillar $P$ using $\{u_1,x_1\}, \{x_1,x_2\} ,\{u_2,x_2\} \in E(G)$ and adding $p$ vertices from $N(x_1) \cap Z - (V(P_1) \cup V(P_2))$ and another $p$ vertices from $N(x_2) \cap Z - (V(P_1) \cup V(P_2))$ as leaf vertices of $x_1,x_2$.
$\Box$

Now, let $G':=G[V\setminus (Z\cup V(P_{abs})]$ and let $P$ be a longest caterpillar in $G'$.
Starting with $P$ we will extend $P$ iteratively, adding at least $\delta C/2$ vertices by using at most $10(p+1)$ vertices from $Z$ in each step, until the number of vertices left is at most $\frac{\delta^2 n}{2}$.
Since the number of iterations is at most $2n/(\delta C)$, and so the number of vertices used to construct $P$ in $Z$ is at most $\frac{2n}{ \delta C} \cdot 10(p+1) < \frac{\beta^3 \gamma}{4} n $, by Claim \ref{gluing}, the process can be completed.
Moreover, $P_{abs}$ can be connected with $P$ using $Z$ and the number of remaining vertices which are not on the caterpillar is at most $|Z| + \frac{\delta^2 n}{2} \leq \delta^2 n$ and so they can be absorbed by $P_{abs}$.
For the general step, let $W:= V(G') \setminus V(P)$ and suppose $|W|> \frac{\delta^2 n}{2}$. We partition $P$ into $l$ blocks $B_1, \dots, B_l$ of consecutive caterpillars so that $C\leq |B_i|\leq (1+\delta)C$.
\begin{claim}\label{cat-claim}
If $|| G[W] || \geq \gamma |W|^2$, then there is a caterpillar in $G[W]$ with at least $\gamma |W| -p$ vertices.
\end{claim}
{\bf Proof.} $G[W]$ contains a subgraph $H$ such that $\delta(H)>\gamma |W|$.
Let $Q$ be a longest caterpillar in $H$. If $|Q|\leq \gamma |W| -p$, then each endpoint of $Q$ has a neighbor $x \in V(H)\setminus Q$, and every vertex not on $Q$ has at least $p$ neighbors outside $Q$. $\Box$

\noindent {\bf Case 1.} $||G[W]||\geq \gamma |W|^2$. \\
By Claim \ref{cat-claim} there is a caterpillar $Q$ in $G[W]$ on at least $\delta C/2$ vertices.
Since $Q \cap Z, P \cap Z = \emptyset$, by Claim \ref{gluing}, we can construct a caterpillar containing both of them.\\
\noindent {\bf Case 2.} There is a block $B_i$ such that $||B_i, W||\geq \left(\frac{1}{2}+\delta\right)|B_i||W|$.\\
Let $W':=\{w\in W| |N(w)\cap B_i| \geq \left(\frac{1}{2} +\frac{\delta}{2}\right)|B_i|\}$.
Then $|W'|\geq \delta |W| \geq \frac{\delta^3n}{2}$. Consider $H:=G[W', B_i]$. Since there are less than $2^{(1+\delta)C}$ subsets of $B_i$ of size $\left(\frac{1}{2}+\frac{\delta}{2}\right)|B_i|$ , there is a set $X \subset B_i$ such that $|X|=\left(\frac{1}{2}+\frac{\delta}{2}\right)|B_i|$ and for at least $|W'|/2^{(1+\delta)C}$ vertices $w\in W'$, $X \subseteq N(w)\cap B_i$.
Since $\frac{|W'|}{2^{(1+\delta)C}} \geq 2C \geq (\frac{1}{2}+ \frac{\delta}{2})|B_i|$, $H$ contains $K_{D,D}$ where $D=\left(\frac{1}{2}+\frac{\delta}{2}\right)|B_i|$ which in turn contains a caterpillar on $2D-p >(\frac{1}{2}+ \frac{\delta}{2})|B_i|$ vertices.
By Claim \ref{gluing}, using at most $4(p+1)$ vertices in $Z$ we can connect the endpoints of this caterpillar with the endpoints of $B_{i-1}$ and $B_{i+1}$.\\
\noindent {\bf Case 3.} For every block $B_i$, $||B_i, W||< \left(\frac{1}{2}+\delta\right)|B_i||W|.$\\
Since we are not in Case 1, $\sum_{v \in W}|N(v) \cap W| < 2\gamma |W|^2$ and so $\sum_{v \in W}|N(v) \cap P| > (1/2  - \delta - 2\gamma)n|W| - 2\gamma |W|^2$, so
$$||P,W||\geq \left(\frac{1}{2}-2\delta\right)n|W|.$$
A block $B$ is called {\it good} if $||B,W|| \geq  \left(\frac{1}{2}-2\sqrt{\delta}\right)|W||B|$.
Let $q$ denote the number of good blocks. We have $q \geq (1-3\sqrt{\delta})\frac{n}{C}$ as otherwise
$$||P,W|| \leq q \left(\frac{1}{2}+\delta\right)(1+\delta)C|W| + (l-q)\left(\frac{1}{2}-2\sqrt{\delta}\right)(1+\delta)C|W|$$
which is less than $\left(\frac{1}{2}-2\delta\right) n|W|$. Using the same argument as in Case 2, for a good block $B_i$ we can find set $C_i\subset B_i$ and $D_i\subseteq W$ such that $G[C_i, D_i] =K_{|C_i|,|D_i|}, |C_i|=\left(\frac{1}{2}-3\sqrt{\delta}\right) C$ and $|D_i|\geq C$. Let $U:=\bigcup (B_i\setminus C_i)$ where the union is taken over good blocks.  We have
$$|U| \geq (1-3\sqrt{\delta})\frac{n}{C} \cdot C - \left(\frac{1}{2} -3\sqrt{\delta}\right) C \frac{n}{C} = \frac{n}{2}.$$
Thus, since $G$ is not $\beta$-extremal, $||G[U]||\geq \beta n^2$, and so there exist two good blocks $B_s, B_t$  with $s < t$ such that $||G[(B_s \setminus C_s)\cup (B_t \setminus C_t)]|| \geq \beta C^2/2$.
Thus by Claim \ref{cat-claim} $G[(B_s \setminus C_s)\cup (B_t\setminus C_t)]$ contains a caterpillar $S$ on $\beta C/4$ vertices.
In addition $G[C_s \cup C_t, W]$ contains two disjoint caterpillars $S_1,S_2$, each on $(1-7\sqrt{\delta})C$ vertices.
Thus $|S \cup S_s \cup S_t| - |B_s \cup B_t| \geq 2(1-7\sqrt{\delta})C + \frac{\beta C}{4} - 2(1+\delta)C \geq \delta C$.
By using at most $10(p+1)$ vertices in $Z$, we can connect the endpoint of $B_{s-1}$ to the endpoint of $S \cap B_s$, connect the endpoint of $S \cap B_t$ to the endpoint of $B_{t-1}$, connect the endpoint of $B_{s+1}$ to the endpoint of $S_1$, and connect the endpoint of $S_1$ to the endpoint of $S_2$.
Finally, by connecting the endpoint of $S_2$ to the endpoint of $B_{t+1}$, we form a longer caterpillar having more than at least $\delta C$ vertices than previous caterpillar. $\Box$

\section{Extremal case}\label{extremal-sec}
In this section we will address the extremal cases. First we will deal the the case when vertices of $G$ can be partitioned into two sets $V_1,V_2$ such that $||V_1,V_2|| \leq \beta n^2$ and so, $G$ is close to a union of two complete graphs and then we address the case when $G$ has a large almost independent set.

We will start with the following lemma.

\begin{lemma}\label{lem-ext-1}
Let $p \in Z^+$. For any $\xi < 1/(4p+5)$ there is $n_0 \in \mathbb{N}$ such that the following holds.
Let $H$ be a graph on $n \geq n_0$ vertices such that $(p+1)| |H|$ and $\delta(H) \geq (1-\xi)n$.
Let $x,y \in V(H)$. Then there is a spanning $p$-caterpillar in $H$ connecting $x$ and $y$.
\end{lemma}
{\bf Proof.}
Let $P$ be a longest $p-$caterpillar in $G$ connecting $x$ and $y$. 
%Let $S=u_1\dots u_q$ denote the spine of $P$ and let $C[i]$ denote the set of spikes of $u_i$. For any two $x,y\in G$, $|N(x)\cap N(y)| \geq (1-2\xi)n$ and so $|S|\geq (1-\xi)n/(2(p+1))$. Indeed, there is an $x,y$-path in $G$ of length  $\lceil (1-\xi)n/(2(p+1)) \rceil$ and each vertex on the path has $p+1$ distinct neighbors outside it, because $(1-\xi)n > p+1 + (1-\xi)n/2$.%
Let $S=(x=)u_1 \dots u_q(=y)$ denote the spine of $P$ and let $C[i]$ denote the set of spikes of $u_i$. 
For any two $u, v \in G$, $|N(u) \cap N(v)| \geq (1-2\xi)n$ and so $q \geq (1-2\xi)n/(p+1)$.
Indeed, if $q < (1-2\xi)n/(p+1)$ then $|V(P)| < (1-2\xi)n$ and then there exists $u_1' \in (N(u_1) \cap N(u_2)) \setminus V(P)$, $|N(u_1') - V(P)| \geq (1-\xi)n - (1-2\xi)n > p$.
%Since $\delta(H) \geq (1-\xi)n$, $|S| \geq \frac{(1-\xi)}{p+1} n$.
If $V(P) = V(H)$ then we are done, so assume that there exists $\{v', y_1,...,y_p\} \subset V(H) \setminus V(P).$
Since $d(v') \geq (1-\xi)n$ there exists $i \in [q]$ such that $u_i,u_{i+1} \in N(v')$. Otherwise, since $\xi < 1/(4p+5)$,
\begin{align*}
(1-\xi )n \leq d(v') \leq (n-q/2) \leq (1- \frac{1-2\xi}{2(p+1)})n \leq (1- \frac{1-\xi}{4(p+1)})n < (1-\xi)n, 
\end{align*}
a contradiction. Moreover,  there are $p$ distinct vertices $f_1,...,f_p \in [q] \setminus \{i,i+1\}$ such that for each $j \in [p]$ , $|N(v') \cap C[f_j]| > 0$ and $f_jy_j\in E(H)$, which gives us the caterpillar $P'$ such that $V(P') = V(P) \cup \{v',y_1,...,v_p\}$ and $P'$ still connects $x$ and $y$.
$\Box$

A $p$-star is a star which has exactly $p$ leaves. 
\begin{lemma}\label{first-extremal-case}
Let $p\in Z^+$. There is $\beta>0$ and $n_0$ such that if $G$ is a graph on $n \geq n_0$ vertices such that $(p+1)|n$, $\delta(G)\geq \frac{n-1}{2}$, and for some partition $V_1, V_2$ of  $V(G)$ with $|V_i|\geq (1/2 -\beta)n$, $||G[V_1,V_2]||\leq \beta n^2$, then $G$ contains a spanning $p$-caterpillar.
\end{lemma}
{\bf Proof.} Let $\xi$ and $\beta$ be such that $0 < \xi < 1/(4p+5)$ and $0 < \beta \leq (\frac{\xi}{30p})^2$.
Let $W_i:=\{v\in V_i : |N(v) \cap V_i|< (1/2 - 5\sqrt{\beta})n \}$. 
We have $\sum_{v\in V_i} |N(v)\cap V_i| \geq (1/4-\beta/2)n(n-1) -\beta n^2 \geq (1/4 - 2\beta)n^2$ and $$\sum_{v\in V_i} |N(v)\cap V_i| < |W_i|(1/2-5\sqrt{\beta})n  + (|V_i|-|W_i|)|V_{i}|$$
and so $|W_i|\leq \sqrt{\beta}n - 1$. In addition, for every $v\in W_i$, $|N(v)\cap (V_{3-i} \setminus W_{3-i})|\geq 4\sqrt{\beta}n$. Let $U_i:=V_i\setminus W_i$ and $X_i:=W_{3-i}$. Then
\begin{itemize} 
\item for every $v \in U_i$, $|N(v) \cap U_i| \geq (1/2-6\sqrt{\beta})n$,
\item for every $v \in X_i$, $|N(v)\cap U_i| \geq 4\sqrt{\beta}n$.
\end{itemize}

Without loss of generality, suppose $|U_1\cup X_1| \leq |U_2 \cup X_2|$. Then for every $v\in U_1\cup X_1$, $|N(v)\cap (U_2\cup X_2)|\geq 1$.
Let $r_i:=|U_i\cup X_i| \bmod (p+1)$. Since very vertex in $U_1\cup X_1$ has at least one neighbor in $U_2\cup X_2$, we pick $r_1$ vertices $u_1, \dots, u_{r_1}$ in $U_1 \cup X_1$ and choose one neighbor in $U_2 \cup X_2$ for each. 
Note that clearly these neighbors do not need to be distinct. Let $w_1, \dots, w_l$ denote distinct vertices in $U_2 \cup X_2$ chosen in this way.
We have $l \leq p$ and each $w_i$ was chosen by at most $r_1\leq p$ vertices.
We will construct a spanning caterpillar by starting with $(U_1 \cup X_1) \setminus \{u_1, \dots, u_{r_1}\}$. Since $|X_1|\leq \sqrt{\beta}n - 1$, there is a matching from $X_1\setminus \{u_1, \dots, u_{r_1}\}$ to $U_1\setminus \{u_1, \dots, u_{r_1}\}$. 
The matching can be easily extended to a caterpillar $P$ in $G[(U_1\cup X_1) \setminus \{u_1, \dots, u_{r_1}\}]$ on at most $2(p+1)\sqrt{\beta}n$ vertices. 
Let $b$ be the starting point of $P$. Let $G' = G[(U_1\cup X_1) \setminus (\{u_1, \dots, u_{r_1}\} \cup V(P)) $ and $b' \in N(b) \cap V(G')$.
Since $\delta(G') \geq (1/2 - (8+2p)\sqrt{\beta})n \geq \frac{(1-\xi)n}{2} \geq (1-\xi)|G'|$, by Lemma \ref{lem-ext-1} $G'$ contains a spanning $p-$caterpillar $P'$ starting at $b'$.

Denote by $x$ the another starting point of $P'$, i.e $P'$ is $b',x$-caterpillar. Now, pick $y \in N(x) \cap (U_2 \cup X_2)$. 
If $y \in X_2$ then construct a star $Y_0$ centered at $y$ such that $Y_0 \subset (X_2 \cup Y_2) \setminus \{w_1, \dots, w_l\}$ and choose $y' \in N(y) \cap (U_2 \setminus (\{w_1, \dots, w_l\} \cup Y_0))$, otherwise $y' = y$.
We will now construct a caterpillar in $G[U_2\cup X_2 \cup \{u_1, \dots, u_{r_1}\}]$. 
If $a_i$ denotes the number of vertices which choose $w_i$, then select $p-a_i$ neighbors of $w_i$ in $U_2 \cup X_2\setminus \{y, w_1, \dots, w_l\}$, all vertices distinct for different values of $i$. 
Let $S_i$ denote the $p$-star with center at $w_i$. Note that $y$ can be among $w_1, \dots, w_l$ but it cannot be among the remaining vertices of $S_1, \dots, S_l$.
Since $|X_2| \leq \sqrt{\beta}n$, there is a matching from $X_2 \setminus ( \{y,y'\} \cup Y_0 \cup \bigcup_{i \in [l]}S_i ) $ to $U_2 \setminus ( \{y,y'\} \cup Y_0 \cup \bigcup_{i \in [l]}S_i )$. 
The matching and $S_1, ..., S_l$ also can be extended to a caterpillar $P''$ in $G[X_2 \cup U_2 \cup \{u_1,\dots, u_{r_1} \} ]$.
Denote by $y''$ the other endpoint  of the spine  of $P''$ and let $y''' \in (N(y'') \setminus (V(P'') \cup Y_0)  ) \cap U_2$. 
Let $G'' = G[U_2 \cup X_2 \setminus (V(P'') \cup \{ w_1, \dots, w_l\})]$. Since $\delta(G'') \geq (1-\xi)|G''|$, again by Lemma \ref{lem-ext-1},
there exists a spanning $p$-caterpillar $P'''$ connecting $y'$ and $y'''$. Then we get a spanning $p$-caterpillar of $G$ by linking $P',P''$ and $P'''$.
$\Box$

We will now proceed to prove the other extremal case. We have the following lemma.
\begin{lemma}\label{lem-ext-2}
Let $p \in Z^+$. For any $\xi < 1/(4p+5)$ there is $n_0 \in \mathbb{N}$ such that the following holds.
Let $H =(A_1,A_2)$ be a bipartite graph on $n \geq n_0$ vertices with $(p+1)|n$ such that  $|A_1|=|A_2|=\frac{n}{2}$ if $n/(p+1)$ is even and $|A_2|=\frac{n+p-1}{2}$ if $n/(p+1)$ is odd.
Suppose that for any $v \in A_i$, $d(v) \geq (1- \xi)|A_{3-i}|$. For any $x \in A_1$, there exists a spanning $p$-caterpillar starting at $x$ in $H$.
\end{lemma}
{\bf Proof.} First suppose $n/(p+1)$ is even. Let $B_i$ be an arbitrary set of $n/(2(p+1))$ vertices in $A_i$ such that $x\in B_1$. For any vertex $v\in B_i$ and any set $C\subseteq A_{3-i}$ of size $n/(2(p+1))$, $|N(v)\cap C| \geq |C| - \xi n >|C|/2$. Consequently, by Hall's theorem, there is a set of pairwise disjoint $p$-stars with centers in $B_i$ and leaves in $A_{3-i}\setminus B_{3-i}$. In addition, $G[B_1,B_2]$ has a Hamilton cycle and so a spanning path which starts at $x$. The path, in connection with stars, gives a $p$-caterpillar starting at $x$.
Now suppose $|A_2|=\frac{n+p-1}{2}$. Let $B_2$ be a subset of $A_2$ of size $(n-p-1)/(2(p+1))$ and let $B_1$ be a subset of $A_1$ of size $(n+p+1)/(2(p+1))$. Note that 
$|B_i|p=|A_{3-i}|-|B_{3-i}|$. As before, by Hall's theorem there are pairwise disjoint $p$-stars with centers in $B_i$ and leaves in   $A_{3-i}\setminus B_{3-i}$ and $G[B_1,B_2]$ has a spanning path.
$\Box$

\begin{lemma}\label{second-extremal-case}
Let $p\in Z^+$. There is $\beta > 0$ and $n_0$ such that if $G$ is a graph on $n\geq n_0$ vertices such that $(p+1)|n$, for some set $S$ of  $V(G)$ with $|S| \geq (1/2 -\beta)n$, $||G[S]||\leq \beta n^2$,  and
$$\delta(G)\geq \left\lbrace \begin{array}{ll} \frac{n}{2}   & \text{ if } n/(p+1) \text{ is even} \\
\frac {n+1}{2} & \text{ if }  n/(p+1) \text{ is odd and } p> 2\\
\frac{n-1}{2} & \text { if } n/(p+1) \text{ is odd and } p\leq 2\end{array}\right.
$$
then $G$ contains a spanning $p$-caterpillar.
\end{lemma}
{\bf Proof.} Let $\xi$ and $\beta$ be such that $0 < \xi < 1/(4p+5)$ and $0 < \beta \leq \min \{(\frac{\xi}{96})^2, (\frac{\xi}{10+3p})^2 \}$. 
We may assume that $|S| \leq n/2$. Let $U_1 := S$ and $U_2:= V \setminus S$. We have
$$||G[U_1, U_2]||\geq (1/2-\beta)n^2/2 -2\beta n^2 \geq (1-10\beta)|U_1||U_2|.$$
Let $W_i:=\{u\in U_i| |N(u) \cap U_{3-i}| \leq (1-10\sqrt{\beta})|U_{3-i}|\}$. Then
$$||G[U_1, U_2]||\leq |W_1||U_2|(1-10\sqrt{\beta}) +(|U_1|-|W_1|)|U_2|$$ 
and so $|W_1|\leq \sqrt{\beta} |U_1|$ and similarly $|W_2|\leq \sqrt{\beta}|U_2|$. 

We define $s$ to be $n/2$ when $n/(p+1)$ is even and $(n-p+1)/2$ when $n/(p+1)$ is odd.  Let $W:=W_1\cup W_2$. Distribute vertices from $W$ to $X_1, X_2$ so that the following holds.
\begin{itemize}
\item[(a)] If $x\in X_i$, then $|N(x) \cap U_{3-i}|\geq 10 \sqrt{\beta}n$.
\item[(b)] $|\min\{|X_1\cup (U_1\setminus W_1)|, |X_2\cup (U_2\setminus W_2)|\}-s|$ is the least possible.
\end{itemize}
If the quantity in the second condition is positive, we further move vertices from $U_i \setminus W_i$ to $X_{3-i}$ which satisfy (a) to make $|\min\{|X_1\cup (U_1\setminus W_1)|, |X_2\cup (U_2\setminus W_2)|\}-s|$ as small as possible. Let $Y_i := X_i \cup (U_i \setminus W_i)$ and suppose $|Y_1|\leq |Y_2|$. 

First, assume that $|Y_1| = s$. Since for each $v \in W_1 \cup W_2$, $d(v) \geq 10\sqrt{\beta}n > |W_1|+|W_2|$, there is a matching $M$ such that for any $e \in M$, $|e \cap W_1| + |e \cap W_2| = 1$.
Then we extend this matching to $p-$caterpillar $P$ so that for any $e \in M$, $e \cap W_i$ is a vertex of spike. Let $G' = G[V \setminus V(P)] = (V',E')$ and note that $V' \cap W = \emptyset$.
Let $x$ be a last vertex of $P$ and $x' \in N(x) \cap V'$. Let $Y_i' = Y_i \cap V'$. For any $v \in Y_i'$,
$$ |N(v) \cap Y'_{3-i}| \geq (1-10\sqrt{\beta})|U_{3-i}| - 3p \sqrt{\beta}|U_{3-i}| \geq (1-\xi)|Y'_{3-i}|.$$
By Lemma \ref{lem-ext-2}, there exists a spanning caterpillar $P'$ starting at $x'$ of $G'$, then we get a spanning caterpillar of $G$ by attaching $P$ to $P'$. 

Now, we assume that $|Y_1| \neq s$. If $|Y_1|>s$ then  $n/(p+1)$ is odd and since $|Y_1|\leq |Y_2|$,  $p \geq 3,$i.e $\delta(G) \geq \frac{n+1}{2}.$ We have $\delta(G[Y_2])\geq \delta(G) - |Y_1| \geq 1$. 
If $|Y_1| < s$ (and so $|Y_1|< n/2$), then $\delta(G[Y_2])\geq \delta(G)-|Y_1| \geq 1$.

In the first case we proceed as follows. Let $l=\frac{n+p-1}{2}-|Y_1|$. Since $|Y_1|> \frac{n-p+1}{2}$, $l< p-1$.
If there is a vertex $y\in Y_2$ such that $|N(y) \cap Y_2|\geq p-1$, then pick $p-l$ neighbors of $y$ from $Y_1$, $l$ from $Y_2$ to form a $p$ star $S$ centered at $y$ and let $x$ be one more neighbor of $y$ in $Y_2$. Deleting $x$ and all vertices in $S$ gives $Y_1', Y_2'$ such that $|Y_i'|=\frac{n-p-1}{2}$, and so by Lemma \ref{lem-ext-2} there is a spanning $p$-caterpillar in $G[Y_1', Y_2']$ starting at $x$.
If no such vertex exists, then $\Delta(G[Y_2])\leq p-2$. Since $\delta(G[Y_2]) \geq 1$, there is a matching in $G[Y_2]$ of size at least $n/2(p-1)(>p+1)$. Let $y \in Y_2$ be arbitrary and let $x$ be a neighbor of $y$ in $Y_2$. 
Let $M=\{a_1b_1, \dots, a_lb_l\}$ be a matching in $G[Y_2]$ such that $x, y \notin V(M)$. We construct caterpillar $Q$ as follows. Start with $x$ and pick $p$ neighbors of $x$ in $Y_1$. We will use $y$ as a vertex on spine of $Q$. 
Pick a neighbor new vertex $y'\in N(y) \cap Y_1$ and a $a_1'\in N(a_1) \cap Y_1$. Note that $\Delta(G[Y_1])\leq 20\sqrt{\beta}n$ as we can't move any vertices from $Y_1$ and so any two vertices in $Y_1$ have at least $n/4$ common neighbors in $Y_2\setminus V(M)$. 
Select one such unused vertex $z$ which gives a $y,a_1$-path of length four which will be added to the spine of $Q$. Now select $p$ neighbors from the opposite set for each vertices except $a_1$.
In the case of $a_1$, pick $p-1$ neighbors from $Y_1$ and make $b_1$ one of the spikes. Now continue to add additional vertices. Then $Q$ has $2l+2$ spine vertices in $Y_2$, $2l$ spine vertices in $Y_1$ and $|V(Q)\cap Y_2|= (2+2l) +2lp +l$, $|V(Q) \cap Y_1| =(2l+2)p +2l -l $. 
This concludes the construction of $Q$. Let $x'$ be one new neighbor of $a_l$ in $Y_1$. Note that
$|Y_2 \setminus V(Q)|= \frac{n-p+1}{2} + l - (2+2l +2lp +l ) = \frac{n'+p-1}{2} $ where $n' = n - (4l+2)(p+1) = |V - V(Q)|$. Thus by Lemma \ref{lem-ext-2} we can extend $Q$ to get a spanning caterpillar in $G$.

In the second case, we have $\delta(G[Y_2]) \geq s-|Y_1|\geq 1$ and because no vertex can be moved from $Y_2$ to $Y_1$,  $\Delta(G[Y_2]) \leq 20 \sqrt{\beta}n$. 
Let $M$ be a maximum matching in $G[Y_2]$ and suppose $|M|< s -|Y_1|$. Then the number of edges in $G[Y_2]$ incident to $V(M)$ is at most  $40\sqrt{\beta}n |M|< 40\sqrt{\beta}n(s-|Y_1|)$, but 
$||G[Y_2]||\geq \frac{|Y_2|}{2} (s-|Y_1|)$, and $|Y_2|\geq 80 \sqrt{\beta}n$. 

The rest of the argument is similar to the previous case.
For every $y \in Y_2$, we have $|N(y) \cap Y_{1}| \geq (1/2 -20\sqrt{\beta})n$. Let $M=\{a_1b_1, \dots, a_qb_q\}$. We move $b_1, \dots, b_q$ from $Y_2$ to $Y_{1}$  so that after moving $|Y_1|=s$. 
Note that $|Y_1|=|Y_2|$ or $|Y_1|= \frac{n-p+1}{2}, |Y_2|=\frac{n+p-1}{2}$. 
Then we extend this matching to  a $p-$caterpillar $P$ so that for any $i \in [q]$, $b_i$ is a  spike in $P$. Let $G' = G[V \setminus V(P)] = (V',E')$.
Let $x$ be the last vertex of $P$ in $Y_2$ and $x' \in N(x) \cap V'$. Let $Y_i' = Y_i \cap V'$. For any $v \in Y_i'$, since $q \leq 4\sqrt{\beta}n$,
$$ |N(v) \cap Y'_{3-i}| \geq (1/2 - 24\sqrt{\beta})n \geq (1-\xi)|Y'_{3-i}|,$$
By Lemma \ref{lem-ext-2}, there exists a spanning caterpillar $P'$ starting at $x'$ of $G'$, then we get a spanning caterpillar of $G$ by attaching $P$ to $P'$.
$\Box$

\end{document}